\documentclass[12pt]{report}
\usepackage{amsfonts,amsmath,amssymb,amsthm}
\usepackage[russian]{babel}
\long\def\comment#1\endcomment{}

\newtheorem{thrm}{ \w{Theorem}{Теорема}}

\newcommand{\lnpi}{\ln\left({\frac{\pi}2}\right)}
\newcommand{\binsum}{\sum\limits_{k=1}^{m}{m\choose k}}
\newcommand{\inftybinsum}{\lim\limits_{m\to\infty}\sum\limits_{k=1}^m{m\choose k}}
\newcommand{\inftysum}[1]{\lim\limits_{#1\to\infty}\sum\limits_{k=1}^{#1}}
\newcommand{\limit}[1]{\lim\limits_{#1\to\infty}}

\newcommand{\w}[2]{#1}
\begin{document}


\centerline{{\bf{\huge \w{New proofs of some formulas of }{Бесконечные произведения и новое доказательство}}}}
\centerline{{\bf{\huge \w{Guillera-Ser-Sondow}{формулы Сера-Сондова для постоянной Эйлера.}}}}
\smallskip
\centerline{{ \w{Vasily Bolbachan}{Василий Болбачан}}
\footnote{ys93@bk.ru; \w{Stavropol, school 6, Russia.}{Ставрополь,
МОУ СОШ 6, Россия.}}}
\bigskip
\small
\w {We present logarithmic series for} {В этой работе мы сформулируем и докажем бесконечные произведения для} $u, \ln u$ \w{and the
Euler-Mascheroni constant}{и постоянной Эйлера-Маскерони}  $\gamma$. 
\w{It was indicated by J. Sondow that Theorem 4 and all proofs are new. All proofs are elementary. We present some conjectures.} {Все
доказательства элементарные. По мнению Дж. Сондова доказательство формулы для u новое. Так же он надеется, что новыми являются теоремы 2 и 3.}

\normalsize

\section*{1. \w{Introduction and main results}{Введение и основные результаты}.}

\begin{thrm} \w{For each real}{Для действительных} $u> 0$
$$1=\inftybinsum\frac{(-1)^{k+1}}{ku+1}.$$
\end{thrm}


{ \bf \w{Remark}{Замечание} 2.} \w{For each real}{Для действительных} $u> 0$
$$u=\inftybinsum\frac{(-1)^{k+1}}k\ln(ku+1)=\sum\limits_{m=1}^{\infty}\binsum
{\frac{(-1)^{k+1}}m}\ln(ku+1).\qquad$$ \w{From the proof of Theorem 1 it follows that for each}{Из
доказательства Теоремы 1, следует, что для любого} $u_0>0$ \w{the convergence in Theorem 1}{Сходимость в Теореме
1} \w{is uniform for}{равномерна для} $u\in[0,u_0]$. \w{Hence the first formula of Remark 2 follows by
integrating the formula of Theorem 1.}{Следовательно первая формула Замечания 2 получается интегрированием
формулы Теоремы 1.} \w{The second formula of Remark 2 is deduced from the first one below.}{Доказательство
второй формулы замечания 2 быдет предоставлено ниже.}
\w{The second formula of Remark 2 is}{Вторая формула замечания 2} [GS06, Theorem 5.3, GS08].
\w{In}{В} [GS06] \w{the proof of Theorem 5.3 is easy: Theorem 5.3 follows from Theorem 5.2 (which is easy) and Example 2.4 (which is easy and
well-known).}{доказательство Теоремы 5.3 довольно просто: оно вытекает из Теоремы 5.2 и примера 2.4(простого и хорошо известного)}

\smallskip {\bf \w{Corollary}{Следствие} 3.} {\it \w{For each real}{Для действительных} $u>0$}
$$u=\lim\limits_{m\to\infty}\prod\limits_{k=1}^m(k+u)^{{m\choose k}(-1)^{k+1}}=\prod\limits_{m=0}^{\infty}\left(\prod\limits_{k=0}^m(k+u+1)^{{m\choose k}(-1)^k}\right).$$

\w{Take}{Положим, последовательно} $u=1, 2, 3$
\w{in the second formula of Corollary 3}{во второй формуле Следствия 3}
$$1=\frac21\cdot\frac23\cdot\frac89\cdot\frac{128}{135}...,\quad
2=\frac31\cdot\frac34\cdot\frac{15}{16}\cdot\frac{125}{128}...,\quad
3=\frac41\cdot\frac45\cdot\frac{24}{25}\cdot\frac{864}{875}...$$

\w{Recall that}{Напомним, что} $\gamma=\lim\limits_{m\to\infty}\left(1+\frac 12+\frac 13+...+\frac 1m-\ln m
\right).$

{\bf  \w{Theorem}{Теорема} 4.}
$$\gamma=\inftybinsum\frac{(-1)^k}k\ln(k!).$$

\w{Using this formula we reprove the following formula}{Используя эту формулу мы также получим}

\smallskip
{\bf  \w{Corollary}{Следствие} 5.}
$$\gamma=
\sum\limits_{j=1}^{\infty}\sum\limits_{i=1}^j
{j-1 \choose i-1}\frac{(-1)^i}j\ln i.$$

\w{The proof of this formula in [Se26] used analytic continuation of Riemann zeta function. 
 Ideas of proofs of this formula in [So03] are explained [So03, remark before Proof 1].}{Доказательство этой формулы в [Se26] аналитическое продолжение Зета функции Римана. Первые два доказательсва этой формулы в [So03] используют гипергеометрическую функцию.} 





\section*{2. \w{Conjecture.}{Проблемы.}}

1. \w{For each real positive numbers}{Для действительных положительных чисел} $z_1,z_2,\dots,z_n$$$z_1z_2\dots z_n=\lim\limits_{m\to\infty}\sum\limits_{k=1}^m{m\choose k}\frac{(-1)^{k+1}}k\ln(1+z_1\ln(1+z_2\dots\ln(1+z_nk)\dots).$$

For example for $n=2$, we obtain

$$z_1z_2=\lim\limits_{m\to\infty}\sum\limits_{k=1}^m{m\choose k}\frac{(-1)^{k+1}}{k}\ln(1+z_1\ln(1+kz_2)).$$

\smallskip
2. \w{For each real}{Для действительных} $z\geq 0$ $$\sum\limits_{n=0}^{\infty}\frac{(-1)^nz^n}{(n!)^2}=\lim\limits_{m\to\infty}\sum\limits_{n=1}^m{m\choose n}\dfrac{(-1)^{n+1}}{z(n+1)+1}.$$

\smallskip
3.

$$\ln\frac\pi2=
\inftybinsum\dfrac{(-1)^k}{k\sum\limits_{n=1}^m\frac{2^{n-1}}n}
\ln\frac{k!!}{(k-1)!!}=
\sum\limits_{n=1}^{\infty}\sum\limits_{k=1}^n{n\choose k}
\frac{(-1)^{k+1}}{2^n}\left(\frac 2k-\frac 2j-\frac jk+2\right)\ln\frac{k!!}{(k-1)!!}.$$

\w{Where}{Где} $k!!=1\cdot 3\cdot 5\cdot...\cdot(k-2)\cdot k$ \w{for}{для нечетного} $k$ \w{odd and}{и} $k!!=2\cdot 4\cdot 6\cdot...\cdot(k-2)\cdot k$, \w{for}{для четного} $k$ \w{even}{}. \w{See also}{Смотрите также} [So05].

\smallskip
4. \w{For all different positive integers}{Для любых различных натуральных чисел} $a_1,a_2,\dots,a_m$, \w{we define}{определим}

$$f(a_k)=\prod\limits_{n=1,n\ne k}^m\dfrac{a_n}{a_n-a_k}.$$

(\w{For}{для} $m=1$ we have $f(a_1)=1$.)
\w{Then}{тогда}
$$\gamma=\sum\limits_{n=1}^mf(a_n)\dfrac{\ln(a_n!)}{a_n}+\sum\limits_{n=1}^{\infty}\int\limits_0^{1/n}\dfrac{dx}{\prod\limits_{k=1}^m1+(a_kx)^{-1}}.$$

\smallskip
5. Define $y_k:=\sum\limits_{n=1}^{\infty}\int\limits_0^{1/n}\frac{dx}{1+(kx)^{-1}}$.
Then the numbers $y_i$, where $i$ runs through positive integers, are linearly independent over $\mathbb{Z}$.

6. It follows from 4 and 5 that $e^{\gamma}$ \w{is a  irrational number}{иррациональное число}.


\smallskip

\section*{3. Proofs.}

\w{All Lemmas are essentially known (Lemma 9 and Lemma 10 can be found in [Wi08]). But we present proofs for
completeness.}{Ни одна из лемм не является принципиально новой (Леммы 9 и 10 можно найти в [Wi08]). Но, для
полноты изложения, мы все же предоставим их доказательства.}

\w{In order to prove Theorem 1 and Theorem 4 we need}{Для доказательства Теоремы 1 и Теоремы 4, нам понадобится}

\smallskip
{\bf \w{Lemma}{Лемма} 7.} {\it \w{For each}{Для действительных} $z>0$ \w{and}{и} $m=1,2,3...$}
$$\sum\limits_{k=0}^m{m\choose k}\frac{(-1)^k}{k+z}=\frac{g_m(z)}z\quad\mbox{\w{where}{где}}\quad
g_m(z):=\frac{m!}{(z+1)(z+2)\dots(z+m)}.$$

{\it \w{Proof}{Доказательство}.}
\w{Take the decomposition}{Возьмем разложение} $\dfrac{g_m(z)}z=\dfrac{A_0}z+\dfrac{A_1}{z+1}+\dfrac{A_2}{z+2}+...+\dfrac{A_m}{z+m}$
\w{into simplest fractions.}{на простейшие дроби}
\w{We have}{Имеем} $m!=\sum\limits_{k=0}^{m}A_k\prod\limits_{i=0,i\ne k}^{m}(z+i)$.
\w{Taking}{положив} $z=-k$ \w{we obtain}{мы получим}
$$m!=A_k(-k)(-k+1)...2\cdot1\cdot2...(m-k-1)(m-k)=A_k k!(m-k)!(-1)^k.$$
\w{Hence}{Слeдовательно} $A_k={m\choose k}(-1)^k$. QED

\smallskip
{\bf  \w{Lemma}{Лемма} 8.} \w{For each}{любых} $z_0>0$ \w{we have that}{выражение} $g_m(z)/z$ \w{converges
to 0 uniformly for}{Сходится к 0 равномерно для} $z\in [z_0;+\infty)$ \w{as $m$ tends to infinity}{при стремлении $m$ к бесконечности}.

\smallskip
{\it \w{Proof}{Доказательство}.} \w{We have}{Имеем}
$$\frac{g_m(z)}z=\frac{1}{z\left(1+z\right)\left(1+\dfrac
z2\right)\dots \left(1+\dfrac zm\right)}< \frac{1}{z\cdot
z\left(1+\dfrac12+\dfrac13+\dots+\dfrac1m\right)}\underset{\
m\to\infty}\to0.\quad QED$$

\smallskip
{\it \w{Proof of Theorem}{Доказательство Теоремы} 1.} \w{It follows from Lemma 7 and Lemma 8
that}{Из Леммы 7 и Леммы 8 следует, что}
$$\limit{m}\sum\limits_{k=0}^m{m\choose k}\frac{(-1)^k}{k+z}=0.$$
\w{Taking}{Положив} $z=\frac 1u$ \w{and changing the limit of summation we obtain}{и поменяв предел суммирования, мы получим}
$$\inftybinsum\frac{(-1)^{k+1}}{ku+1}=1.\quad QED$$

\smallskip
{\it \w{Proof of Theorem 4.}{Доказательство Теоремы 4.}}
\w{Let us prove that}{Докажем, что}
$$\binsum\frac{(-1)^{k+1}}{k}\ln(k!)\overset{(1)}=\lim\limits_{n \to \infty}\left(\sum\limits_{k=1}^m{m\choose k}(-1)^{k+1}\ln n-\sum\limits_{j=1}^n\binsum\frac{(-1)^{k+1}}{k}\ln\left(1+\frac kj\right) \right)\overset{(2)}=$$
$$\overset{(2)}=\lim\limits_{n\to\infty}\left(\ln n-\sum\limits_{j=1}^n\frac 1j\right)+\sum\limits_{j=1}^{\infty}\int\limits_0^{1/j}g_m(1/u)du\overset{(3)}=-\gamma+\sum\limits_{j=1}^{\infty} \int\limits_{0}^{1/j}
g_m(1/u)du \overset{(*)}\to-\gamma\quad\mbox{as}\quad m\to\infty.$$
\w{The first equality follows because}{Первое равенство верно т.к.}
$$k!=k!\lim\limits_{n\to\infty}\frac{n^k}{(n+1)(n+2)\dots(n+k)}=\limit{n}\frac{n^kn!k!}{(n+k)!}=$$$$=\limit{n}\frac{n^kn!}{(k+1)(k+2)\dots(k+n)}=\limit{n}\frac{n^k}{\left(1+\frac k1\right)\left(1+\frac k2\right)\dots\left(1+\frac{k}n\right)}.$$
\w{Let us prove the second equality. Taking $z=\frac 1u$ in Lemma 7 and changing the limit of summation we have}{Докажем второе равенство. Положив $z=\frac 1u$ в Лемме 7 и поменяв предел суммирования, имеем}
$$1-\sum\limits_{k=1}^m{m\choose k}\frac{(-1)^{k+1}}{ku+1}=g_m(1/u).$$
\w{Hence (improperly) integrating this formula with from $0$ to $1/j$, we obtain}{Следовательно (несобственно) интегрируя эту формулу от $0$ до $1/j$, мы получим}
$$\binsum\frac{(-1)^{k+1}}k\ln(1+\frac kj)=\frac 1j-\int\limits_{0}^{1/j}g_m(1/u) du.$$
\w{This and}{Это и то, что} $\sum\limits_{k=1}^m{m\choose k}(-1)^{k+1}=1$ \w{imply the second equality}{доказывает второе равенство}.

\w{The third equality is clear.}{Третье равенство очевидно}

\w{Let us prove (*). We have}{Докажем (*). Имеем}
$$0<g_m(1/u)<\dfrac{g_{m-1}(1/u)}{1+\frac 1m}\quad\mbox{for}\quad 0<u\leq 1.$$
\w{Hence}{Следовательно}
$$0<g_m(1/u)\leq\dfrac{g_1(1/u)}{(1+\frac 12)(1+\frac 13)...(1+\frac 1m)}=\frac{2g_1(1/u)}{m+1}\quad\mbox{for}\quad 0<u\leq 1.$$

\w{For each $m$ the series of left-hand side of the equality (*) converges by the sum of limits theorem.  Hence for the sum $S_m$ in the left-hand side of (*) we have}{Для любых $m$ ряд в левой части равенства (*) Сходится, по теореме о сумме пределов. Следовательно для сумм $S_m$ в левой части (*), имеем}

$$0<S_m<\frac{2S_1}{m+1}\overset{m\to\infty}\to 0.\quad QED$$



\smallskip
\w{In order to prove Remark 2, Corollary 3 and Corollary 5 we need}{Для доказательства Замечания 2, Следствия 3 и Следствие 5 нам понадобится}

\smallskip
{\bf \w{Lemma}{Лемма} 9.}
$$\frac{{m \choose k}}k=\sum\limits_{n=k}^{m}\frac{{n\choose k}}n.$$


\smallskip
{\it \w{Proof.}{Доказательство.}} \w{We have}{Имеем}
$${m \choose k}=\sum\limits_{n=1}^{m-k+1}{m-n\choose k-1}
=\sum\limits_{n=k}^{m}{n-1\choose k-1}=k\sum\limits_{n=k}^{m}\frac{{n\choose k}}n.$$
\w{Here the first equality holds because}{Здесь первое равенство верно, так как } ${m-n}\choose {k-1}$ \w{equals
to the number of $k$-subsets of}{равно числу k - подмножеств} $\{1,2,\dots,m\}$ \w{whose minimal element is $n$,
the second equality holds because the summands in those sums are equal.}{чей минимальный элемент есть n, второе равенство верно, так как слогаемые в обоих суммах одни и те же.}

\smallskip
{\it \w{Proof the second formula of Remark 2}{Доказательство второй формулы Замечания 2}.} \w{By Lemma 9 for}{Используя Лемму 9, для}
$X_{n,k}={n\choose k }\frac{(-1)^{k+1}}n\ln(ku+1)$ \w{we have}{имеем}
$$\inftybinsum\frac{(-1)^{k+1}}k\ln(ku+1)
=\lim\limits_{m\to\infty}\sum\limits_{k=1}^m\sum\limits_{n=k}^mX_{n,k}=\sum\limits_{n=1}^{\infty}\sum\limits_{k=1}^nX_{n,k}.\quad QED$$

\smallskip
\smallskip
{\it \w{Proof the first formula of Corollary 3.}{Доказательство первой формулы Следствия 3.}}
\w{Take logarithms of both sides. Take $z=\frac 1u$ in the first formula of Remark 2 we obtain }{Докажем логарифмический вариант. По Теореме 1}
$$\frac 1z=\limit{m}\sum\limits_{k=1}^m{m\choose k}\frac{(-1)^{k+1}}k\ln\left(\frac kz+1\right).$$

\w{By Lemma 8 the limit in this formula is uniform for}{По Лемме 8 предел (в Теореме 1) сходятся равномерно для} $z\in [z_0;+\infty], z_0\geq 0$. \w{Hence integrating this formula from 1 to u with respect to $z$, we get }{Следовательно интегрируя это равенство по $u$, получим}
$$\ln u\overset{(1)}=\limit{m}\sum\limits_{k=1}^m{m\choose k}\frac{(-1)^{k+1}}k\left(k\ln(k+u)-k+u\ln\left(\frac ku+1\right)\right)\overset{(2)}=$$$$\overset{(2)}=\limit{m}\sum\limits_{k=1}^m{m\choose k}(-1)^{k+1}\ln(k+u)-\limit{m}\sum\limits_{k=1}^m{m\choose k}(-1)^{k+1}+u\limit{m}\sum\limits_{k=1}^m{m\choose k}\frac{(-1)^{k+1}}k\ln\left(\frac ku+1\right)\overset{(3)}=$$$$\overset{(3)}=\limit{m}\sum\limits_{k=1}^m{m\choose k}(-1)^{k+1}\ln(k+u).$$

The first equality follows because

$$\int\limits_1^{u}\frac {dz}z=\ln u,\quad \int\limits_{1}^u\ln(\frac kz+1)dz=(k+u)\ln(k+u)-k+u\ln u=k\ln(k+u)-k+u\ln\left(\frac ku+1\right)$$

The second equality equality is clear. The third equality follows by the second formula of Remark 2 and
$$\sum\limits_{k=1}^{m}{m\choose k}(-1)^{k+1}=1. \quad QED.$$
\smallskip
{\it \w{Proof the second formula of Corollary 3.}{Доказательство второй формулы следствия 3.}} \w{Take logarithms of both sides.}{Докажем логарифмический варинт.} \w{By Lemma 9 for}{Используя Лемму 9, для}
$a_k=\ln(k+u)$ \w{we obtain}{мы получим}
$$\inftybinsum(-1)^{k+1}a_k=\inftysum{m}\sum\limits_{n=k}^m{n\choose k }(-1)^{k+1}\frac{k}na_k=$$$$=\sum\limits_{1\leq k\leq n<\infty}{n-1\choose k-1 }(-1)^{k+1}a_k=\sum\limits_{n=1}^{\infty}\sum\limits_{k=1}^n{n-1\choose k-1 }(-1)^{k+1}a_k=\sum\limits_{n=0}^{\infty}\sum\limits_{k=0}^n{n\choose k }(-1)^ka_{k+1}.\quad QED$$


\smallskip

\w{In order to prove Corollary 5, we need}{Для доказательства Следствия 5, нам понадобится}

\smallskip
{\bf \w{Lemma}{Лемма} 10.}
$$\sum\limits_{k=n}^j(-1)^{k+n}{j\choose k}={j-1\choose n-1}.$$


{\it \w{Proof.}{Доказательство.}}
\w{The proof is by induction on $n$. For $n$ = 1 this is known.}{Докажем индукцией по $n$. База индукции при $n=1$ очевидна.}

\w{Let us prove the inductive step. By Pascal's rule we have:}{Докажем шаг индукции. По правилу Паскаля, имеем:}
$$\sum\limits_{k=n+1}^j{j\choose k}(-1)^{k+n}=\sum\limits_{k=n}^j{j\choose k}(-1)^{k+n}-{j\choose n}={j-1\choose n-1}-{j\choose n}={j-1 \choose n}.\quad QED$$

\smallskip
{\it \w{Proof of Corollary}{Доказательство следствия} 5.} We have
$$\gamma\overset{(1)}=\inftybinsum\frac{(-1)^k}k\ln(k!)
\overset{(2)}=\inftysum{m}(-1)^k\left(\sum\limits_{j=k}^m{j\choose
k}\frac{1}j\right)\left(\sum\limits_{n=1}^k\ln n \right)\overset{(3)}=$$
$$\overset{(3)}=\lim\limits_{m\to\infty}\sum\limits_{k=1}^m\sum\limits_{j=k}^m\sum\limits_{n=1}^ka_{kjn} \overset{(4)}=\limit{m}\sum\limits_{1 \leq n\leq k \leq j\leq m}a_{kjn} \overset{(5)}=\limit{m}\sum\limits_{j=1}^m\sum\limits_{n=1}^j\sum\limits_{k=n}^ja_{kjn} \overset{(6)}=$$
$$\overset{(6)}=\sum\limits_{j=1}^{\infty}\sum\limits_{n=1}^j\frac{(-1)^n}j{j-1 \choose n-1}\ln n .$$
\w{Here}{Здесь} $a_{kjn}={j\choose k}\frac{(-1)^k}j\ln n$.
\w{The first equality follows by theorem 4. The second equality
follows by lemma 9. The third, the fourth and the fifth equality is clear. The sixth
equality follows by lemma 10.}{Первое равенство следует из Теоремы 4. Второе равенство следует из Леммы 9. Третье, Четвертое и пятое равенство очевидны. Шестое равенство следует из Леммы 10.} QED

\smallskip
\section*{4. \w{Acknowledgements.}{Благодарности}}
\w{The author is grateful to J. Sondow and A.Skopenkov for useful comments
on this paper and references.}{Автор Блогадарен Дж. Сондову, А.Скопенкову и В.Зудилину за многочисленные замечания по этой работе и ссылкам.}
\w{This paper was presented at Moscow Mathematical Conference of High-school Students in 2009.}{Эта работа была представлена на семинарах А.Скопенкова, А.Сгибнева и на Московской математической конференции школьников в 2009 году.}

\bigskip
\centerline{\bf \w{References}{Ссылки}}

[GS06] J. Guillera and J. Sondow, Double integrals and infinite products for
some classical constants via analytic continuations of Lerch's transcendent,
Ramanujan Journal 16 (2008) 247-270.
http://arxiv.org/abs/math.NT/0506319.

[GS08] J. Guillera and J Sondow, Problem 11381, Amer. Math. Monthly 115 (2008)
665.


[Se26] J.Ser, Sur une expression de la fonction $\zeta\, (s)$ de Riemann, Comptes Rendus 182 (1926) 1075-1077.

[So05] J. Sondow, A faster product for $\pi$ and a new integral for $\lnpi$, Amer. Math. Monthly 112 (2005) 729-734.
http://arxiv.org/abs/math.NT/0401406

[So03]J. Sondow, An Infinite Product for $e^\gamma$ via
Hypergeometric Formulas for Euler's Constant $\gamma$.
http://arxiv.org/abs/math.NT/0306008

[Wi08] http://en.wikipedia.org/wiki/Binomial$\_$coefficient
 \end{document}